\documentclass[12pt,a4paper]{article}
\usepackage{amsmath,amssymb,amsthm}
\usepackage[margin=2.5cm]{geometry}
\usepackage[hidelinks]{hyperref}

\numberwithin{equation}{section}
\newtheorem{theorem}{Theorem}[section]
\newtheorem{lemma}[theorem]{Lemma}
\newtheorem{proposition}[theorem]{Proposition}
\newtheorem{corollary}[theorem]{Corollary}
\theoremstyle{definition}
\newtheorem{definition}[theorem]{Definition}
\newtheorem{remark}[theorem]{Remark}
\newtheorem{result}[theorem]{Result}

\newcommand{\Z}{\mathbb{Z}}
\newcommand{\Q}{\mathbb{Q}}
\newcommand{\C}{\mathbb{C}}

\newcommand{\F}{\mathbb{F}}
\newcommand{\PP}{\mathfrak{P}}
\newcommand{\BH}{\mathrm{BH}}
\newcommand{\ov}{\overline}
\newcommand{\fh}{\widehat{f}}
\DeclareMathOperator{\ord}{ord}

\begin{document}

\title{Classification of circulant Butson matrices of order $p^2$}
\author{Ka Hin Leung\thanks{Department of Mathematics, National University of Singapore, Kent Ridge,
Singapore 119260, Republic of Singapore. \texttt{matlkh@nus.edu.sg}. Research partially supported by
National Science Foundation of China grant 12131011.}
\and Bernhard Schmidt\thanks{Division of Mathematical Sciences, School of Physical and Mathematical
Sciences, Nanyang Technological University, Singapore 637371, Republic of Singapore.}
\and Meng Zhao\thanks{College of General Education, Guangdong University of Science and Technology,
Dongguan, Guangdong, China 523083.}}
\date{}
\maketitle

\begin{abstract}
For every prime $p$, we determine all Butson Hadamard matrices of order $p^2$ invariant under $\Z_{p^2}$.
They are given by an explicit group ring formula and, in particular, 
there are exactly $p!\,h^p$ such matrices with entries in the $h$-th roots of unity if
$p$ divides $h$, and none otherwise.
The proof combines the following ingredients: 
(i) an application of the polynomial method over a local ring $\F[\Z_{p^2}]$, where $\F$ is a field of characteristic $p$,
(ii) the field descent method,
(iii) a combinatorial analysis of the blocks of the matrix based on the geometry of the affine plane ${\rm AG}(2,p)$, and
(iv) Tao's uncertainty principle for cyclic groups of prime order.
\end{abstract}

\medskip
\noindent\textbf{Keywords:}  complex Hadamard matrices, circulant matrices, polynomial method, field descent, uncertainty principle

\noindent\textbf{MSC 2020:} 05B20, 11R18, 05B10.


\section{Introduction}\label{sec:intro}

Throughout, $\zeta_n=e^{2\pi i/n}$ and $\mu_n=\langle\zeta_n\rangle$ is the group of complex
$n$-th roots of unity. An $n\times n$ matrix $M$ whose entries are complex $h$-th roots of unity is
a \emph{Butson Hadamard matrix} if $MM^*=nI_n$, where $M^*$ is the conjugate transpose of $M$.
If $G$ is a finite group, a matrix
$M=(m_{x,y})_{x,y\in G}$ indexed by $G$ is \emph{$G$-invariant} if $m_{xz,yz}=m_{x,y}$ for all
$x,y,z\in G$. A $G$-invariant Butson Hadamard matrix with entries in $\mu_h$ is called a
\emph{$\BH(G,h)$ matrix}. The set of all $\BH(G,h)$  matrices is also denoted by $\BH(G,h)$ .
For a cyclic group $G$ the $G$-invariant matrices are circulant.

For an abelian group $G$ it is convenient to translate this into the language of group rings, as
in \cite{DS}. For a commutative ring $R$ and a finite abelian group $G$, written
multiplicatively with identity $e$, the group ring $R[G]$ consists of the formal sums
$X=\sum_{x\in G}c_xx$ with $c_x\in R$; for $R\subseteq\C$ we put
$X^{(-1)}=\sum_{x\in G}\ov{c_x}\,x^{-1}$. Every character $\chi$ of $G$ is extended linearly to
$\C[G]$. A $G$-invariant matrix has the form $M=(a_{xy^{-1}})_{x,y\in G}$, and by
\cite[Lemma~2.1]{Sch19}, if $a_x\in\mu_h$ for all $x$, then $M$ is a $\BH(G,h)$ matrix if and only
if $D=\sum_{x\in G}a_xx\in\Z[\zeta_h][G]$ satisfies
\[
DD^{(-1)}=|G| ,
\]
if and only if $|\chi(D)|^2=|G|$ for every character $\chi$ of $G$. 
We identify $M$ with $D$ and write $D\in\BH(G,h)$.
Two group ring elements $D,E$ are \emph{equivalent} if
$E=\pm\zeta_h^i\,x\,\alpha(D)$ for some integer $i$, some $x\in G$ and some automorphism $\alpha$ of
$G$.

For a subgroup $U$ of a finite abelian group, we also write $U$ for the group ring element
$\sum_{u\in U}u$. A prime $p$ is \emph{semiprimitive modulo $h$} if
$p^j\equiv-1\pmod h$ for some positive integer $j$.

Recent work has led to substantial progress in the classification of group invariant
Butson Hadamard matrices:
\begin{result}\label{res:known}
Let $p$ be a prime and $h$ a positive integer.
\begin{itemize}
\item[(a)] \cite{GS,LCZ} Write $\Z_p=\langle g\rangle$.
There is a $\BH(\Z_p,h)$ matrix if and only
if $p\mid h$, and $4\mid h$ in case $p=2$. If $p$ is odd and $p\mid h$, then
$D\in\BH(\Z_p,h)$ if and only if
\[
D=\omega\sum_{y=0}^{p-1}\zeta_p^{ay^2+by}g^y
\]
with $\omega\in\mu_h$ and $a,b\in\Z_p$, $a\ne0$. In particular $|\BH(\Z_p,h)|=p(p-1)h$.

\item[(b)] \cite{DS26,SWX21} Let $p$ be odd and let $U_0,\dots,U_p$ be the subgroups of $\Z_p\times\Z_p$
of order $p$. If $\eta_0,\dots,\eta_p\in\mu_h$ satisfy $\sum_{i=0}^p\eta_i=1$, then
\[
D=\sum_{i=0}^{p}\eta_iU_i\in\BH(\Z_p\times\Z_p,h).
\]
If $p$ is semiprimitive modulo $h$, then every $\BH(\Z_p\times\Z_p,h)$ matrix is equivalent to one of this form. 
\end{itemize}
\end{result}

The classification of $\BH(\Z_{p^2},h)$ matrices, however, turns out to be more difficult 
and requires a more advanced toolset.  The following is the main result of this paper.

\begin{theorem}\label{thm:main}
Let $p$ be a prime, let $h$ be a positive integer, and write $\Z_{p^2}=\langle g\rangle$.
A $\BH(\Z_{p^2},h)$ matrix exists if and only if $p\mid h$. When $p\mid h$, the
$\BH(\Z_{p^2},h)$ matrices are precisely the elements
\begin{equation}\label{eq:family}
D=\sum_{t=0}^{p-1}\omega_t\,g^t\sum_{y=0}^{p-1}\zeta_p^{\pi(t)y}g^{py}
\end{equation}
with $\omega_0,\dots,\omega_{p-1}\in\mu_h$ and a permutation $\pi$ of $\Z_p$.
In particular, $|\BH(\Z_{p^2},h)|=p!\,h^p$ if $p\mid h$.
\end{theorem}

\begin{remark}
The number $p!\,h^p$ counts all group ring elements in $\BH(\Z_{p^2},h)$, not equivalence classes.
\end{remark}

For the rest of the paper, $p$ denotes a prime, $G=\langle g\rangle$ is cyclic of order $p^2$,
$w=g^p$, and $H=\langle w\rangle$ is the subgroup of $G$ of order $p$.

\paragraph{Outline of the proof.} The proof consists of the following steps.

\emph{(1) The local-ring theorem} (Section~\ref{sec:local}). Write every entry as
$a_x=\zeta_{p^c}^{a(x)}\eta(x)$ with $\eta(x)\in\mu_{h'}$, where $h=p^ch'$ and $p\nmid h'$.
Reducing the identity $DD^{(-1)}=p^2$ modulo a prime ideal above $p$ gives a relation in the local
ring $\F[G]\cong\F[\rho]/(\rho^{p^2})$, where $\F$ is a finite field of characteristic $p$, and a computation
of orders of elements of this ring shows that $\eta$ is constant on every coset of $H$. This holds
for every prime $p$ and every $h$, and implies $p\mid h$ (Corollary~\ref{cor:pdividesh}).

\emph{(2) Field descent} (Section~\ref{sec:descent}). For odd $p$, the field descent 
\cite{Sch99} shows that, for $c\ge2$, each character value of a $\BH(G,p^ch')$ matrix lies in a coset
$\zeta_{p^c}^e\,\Q(\zeta_{ph'})$, with $e$ depending on the character. For $c=2$ this means that exactly one component of a character
value with respect to the basis $1,\zeta_{p^2},\dots,\zeta_{p^2}^{p-1}$ of $\Q(\zeta_{p^2h'})$ over
$\Q(\zeta_{ph'})$ is nonzero (Lemma~\ref{cor:unique}).

\emph{(3) The case $\nu_p(h)=1$ for odd $p$} reduces to the case $\nu_p(h)=2$, since
$\BH(G,h)\subseteq\BH(G,ph)$.

\emph{(4) The case $\nu_p(h)=2$ for odd $p$} (Section~\ref{sec:nu2}) is the heart of the paper. We
split $D$ into blocks according to the coset of $H$ and the residue modulo $p$ of the exponent
$a(x)$, and regard the blocks as points of the affine plane $\Z_p^2$. For every line,
Lemma~\ref{cor:unique} gives an identity in the group ring of $H$ (Lemma~\ref{lem:lineid}), which,
together with the uncertainty principle for $\Z_p$, bounds the number of elements on the line from
below (Lemma~\ref{lem:line}). According to the minimal number
$A$ of elements of a nonzero block, there are three cases: $A=1$, or $2\le A\le(p-1)/2$, or $A=p$.
The first two are ruled out by counting arguments; in the third, every coset of $H$ carries a single
block, and a count of the directions determined by the exponents of the blocks
(Lemma~\ref{lem:directions}) shows that every block is a constant times a character of $H$.

\emph{(5) The case $\nu_p(h)\ge3$ for odd $p$} (Section~\ref{sec:nu3}) is reduced to $\nu_p(h)=2$ by the field
descent: the entries of $D$ are grouped according to their $\zeta_{p^c}$-exponent modulo
$p^{c-2}$, and the resulting pieces are recombined with a root of unity of a suitable prime order
$q$ in place of $\zeta_{p^c}$.

\emph{(6) The prime $2$} (Section~\ref{sec:p2}). An elementary argument for arbitrary unimodular
coefficients settles this case for every even $h$, including $\nu_2(h)=1$.

The proof of Theorem~\ref{thm:main} is assembled in Section~\ref{sec:proofmain}, where we also
verify the converse and the count. The divisibility $p\mid h$, previously established in
\cite[Theorem~3.6]{DS}, is also proved here as a consequence of the local-ring theorem.
The only results used without proof in the classification argument are the field descent theorem
\cite[Theorem~3.5]{Sch99}, the uncertainty
principle for cyclic groups of prime order \cite{Tao05}, Fourier inversion and Parseval's identity on
finite abelian groups, and standard facts on the arithmetic of cyclotomic fields.

\section{Preliminaries}\label{sec:prelim}

\subsection{Characters, Fourier inversion and Parseval's identity}

Let $K$ be a finite abelian group and $\widehat K$ its character group. For
$X=\sum_{x\in K}c_xx\in\C[K]$ the \emph{Fourier inversion formula}
$c_x=|K|^{-1}\sum_{\chi\in\widehat K}\chi(Xx^{-1})$ holds; in particular, $X$ is determined by
the values $\chi(X)$ with $\chi\in\widehat K$. Since $\chi(X^{(-1)})=\ov{\chi(X)}$, an element
$D=\sum_xa_xx$ with $a_x\in\mu_h$ lies in $\BH(K,h)$ if and only if $|\chi(D)|^2=|K|$ for all
$\chi\in\widehat K$. \emph{Parseval's identity} states that
\[
\sum_{\chi\in\widehat K}|\chi(X)|^2=|K|\sum_{x\in K}|c_x|^2 .
\]
We also use the uncertainty principle for cyclic groups of prime order, due to Tao \cite{Tao05}.

\begin{lemma}\label{lem:uncertainty}
Let $K$ be cyclic of prime order $p$ and let $X=\sum_{x\in K}c_xx\in\C[K]$ be nonzero. Then
\[
|\{x\in K:\ c_x\ne0\}|+|\{\chi\in\widehat K:\ \chi(X)\ne0\}|\ \ge\ p+1 .
\]
\end{lemma}

\subsection{Cyclotomic fields}

\begin{lemma}\label{lem:basis}
Let $p$ be a prime and $h'\ge1$ with $p\nmid h'$.
\begin{itemize}
\item[(a)] For $n\ge r\ge1$, the elements $\zeta_{p^n}^e$ with $0\le e<p^{n-r}$ form a basis of
$\Q(\zeta_{p^nh'})$ over $\Q(\zeta_{p^rh'})$.
\item[(b)] Let $b_{s,k}\in\Q(\zeta_{h'})$ for $s,k\in\{0,\dots,p-1\}$. Then
$\sum_{s,k}\zeta_{p^2}^{s+pk}b_{s,k}=0$ if and only if $b_{s,k}$ does not depend on $k$ for every
$s$. In particular, $\sum_{k=0}^{p-1}\zeta_p^kb_k=0$ with $b_k\in\Q(\zeta_{h'})$ if and only if all
$b_k$ are equal.
\end{itemize}
\end{lemma}

\begin{proof}
(a) We have $\Q(\zeta_{p^nh'})=\Q(\zeta_{p^rh'})(\zeta_{p^n})$ and
$[\Q(\zeta_{p^nh'}):\Q(\zeta_{p^rh'})]=\varphi(p^nh')/\varphi(p^rh')=p^{n-r}$. As $\zeta_{p^n}$ is a
root of $X^{p^{n-r}}-\zeta_{p^r}\in\Q(\zeta_{p^rh'})[X]$, this polynomial is its minimal polynomial,
and the claim follows.
(b) As $\gcd(p,h')=1$, we have $[\Q(\zeta_{p^2h'}):\Q(\zeta_{h'})]=p(p-1)$, so the minimal polynomial of
$\zeta_{p^2}$ over $\Q(\zeta_{h'})$ is $\Phi_{p^2}(X)=\sum_{k=0}^{p-1}X^{pk}$. Hence the
$\zeta_{p^2}^{s+pk}$ with $0\le s\le p-1$ and $0\le k\le p-2$ form a basis, and
$\zeta_{p^2}^{s+p(p-1)}=-\sum_{k\le p-2}\zeta_{p^2}^{s+pk}$. Thus the sum equals
$\sum_s\sum_{k\le p-2}\zeta_{p^2}^{s+pk}(b_{s,k}-b_{s,p-1})$, which vanishes if and only if
$b_{s,k}=b_{s,p-1}$ for all $s,k$. The last claim is the same argument for $\zeta_p$, whose
minimal polynomial over $\Q(\zeta_{h'})$ is $\sum_{k=0}^{p-1}X^k$.
\end{proof}

\begin{lemma}\label{lem:sums}
Let $n_0,\dots,n_{p-1}\in\Z$. Then $\sum_kn_k\zeta_p^k=0$ if and only if $n_0=n_1=\dots=n_{p-1}$.
Consequently, a sum of $n$ (not necessarily distinct) $p$-th roots of unity with $1\le n\le p-1$ is
nonzero.
\end{lemma}

\begin{proof}
The first claim is Lemma~\ref{lem:basis}(b) with $h'=1$. A sum of $n$ $p$-th roots of unity is
$\sum_kn_k\zeta_p^k$ with $n_k\ge0$ and $\sum_kn_k=n$; if it vanishes, then all $n_k$ equal some
$c$, and $n=pc$ is divisible by $p$.
\end{proof}

We will need the following special case of the field descent method.

\begin{result}\label{res:descent}
Let $p$ be an odd prime, let $n,b,h'$ be positive integers with $p\nmid h'$, and let
$X\in\Z[\zeta_{p^nh'}]$ satisfy $X\ov X=p^b$. Then
\begin{equation}\label{eq:coset}
X\in\zeta_{p^n}^{e}\,\Q(\zeta_{ph'})\qquad\text{for some }e\in\{0,\dots,p^n-1\}.
\end{equation}
\end{result}

This is a consequence of the field descent theorem \cite[Theorem~3.5]{Sch99}; see also
\cite[Result~3.4]{DS} for the special case needed here.

\subsection{Polynomial functions and local group rings}

Let $\F$ be a field of characteristic $p$. Every function $\gamma\colon\F_p\to\F$ is given by a
unique polynomial $P_\gamma(Y)=\sum_{k=0}^{p-1}P_kY^k\in\F[Y]$ of degree at most $p-1$ (Lagrange
interpolation). Explicitly,
\[
P_\gamma(Y)=\sum_{k=0}^{p-1}\gamma(k)\bigl(1-(Y-k)^{p-1}\bigr).
\]
We put $\deg\gamma:=\deg P_\gamma$ (with $\deg0=-\infty$).
Note that $P_\gamma(k)=\gamma(k)$ for $k=0,\ldots,p-1$.
We identify $\Z_p$
with $\F_p$ and with $\{0,\dots,p-1\}$ when convenient, and use the convention $0^0=1$.

\begin{lemma}\label{lem:powersum}
\begin{itemize}
\item[(a)] For $n\ge0$: $\sum_{y\in\F_p}y^n=-1$ if $n>0$ and $(p-1)\mid n$, and
$\sum_{y\in\F_p}y^n=0$ otherwise.
\item[(b)] For $\gamma\colon\F_p\to\F$ and $0\le a\le p-2$: $\sum_{y}\gamma(y)y^a=-P_{p-1-a}$.
Hence, for $1\le R\le p-1$, we have $\sum_y\gamma(y)y^a=0$ for all $0\le a<R$ if and only if
$\deg\gamma\le p-1-R$.
\item[(c)] If $\gamma_1\gamma_2\equiv1$ on $\F_p$ and $\deg\gamma_1+\deg\gamma_2\le p-1$, then
$\gamma_1$ and $\gamma_2$ are constant.
\end{itemize}
\end{lemma}

\begin{proof}
(a) For $n=0$ the sum is $p\cdot1=0$. For $n>0$ it equals $\sum_{y\in\F_p^*}y^n$, which is
$p-1=-1$ if $(p-1)\mid n$; otherwise $c^n\ne1$ for a generator $c$ of $\F_p^*$, and the sum $S$
satisfies $c^nS=S$, so $S=0$.
(b) We have $\sum_y\gamma(y)y^a=\sum_kP_k\sum_yy^{a+k}$ and $0\le a+k\le2p-3$; by (a) only $a+k=p-1$
contributes, with $\sum_yy^{p-1}=-1$. The second statement follows, as $a<R$ runs through the
indices $p-1-a>p-1-R$.
(c) The polynomial $P_{\gamma_1}(Y)P_{\gamma_2}(Y)-1$ has degree at most $p-1$ and $p$ roots, so
$P_{\gamma_1}(Y)P_{\gamma_2}(Y)=1$ in $\F[Y]$; hence $P_{\gamma_1}(Y)$ and $P_{\gamma_2}(Y)$ are nonzero
constants in $\F[Y]$. Thus, both $\gamma_1$ and $\gamma_2$ are constant.
\end{proof}

\begin{lemma}\label{lem:order}
\begin{itemize}
\item[(a)] Let $C=\langle c\rangle$ be cyclic of order $p^n$ and $\rho=c-1\in\F[C]$. Then
$\F[C]=\F[\rho]/(\rho^{p^n})$. Every nonzero $\alpha\in\F[C]$ can be written as $\alpha=\rho^t u$
with a unit $u$ and a unique $t=\ord(\alpha)\in\{0,\dots,p^n-1\}$, namely the least exponent
occurring in the representation of $\alpha$ as a polynomial in $\rho$ of degree less than $p^n$.
For nonzero $\alpha,\beta\in \F[C]$ we have $\alpha\beta=0$ if and only if $\ord\alpha+\ord\beta\ge p^n$;
otherwise $\ord(\alpha\beta)=\ord\alpha+\ord\beta$.
\item[(b)] Let $\gamma\colon\Z_p\to\F$ be nonzero, and let
$\sigma_l=\sum_{t=0}^{p-1}\gamma(t)\binom tl\in\F$ for $0\le l\le p-1$. Then
$\min\{l:\sigma_l\ne0\}=p-1-\deg\gamma$. Equivalently, in $\F[C]$ with $n=1$ we have
$\ord\bigl(\sum_{t=0}^{p-1}\gamma(t)c^t\bigr)=p-1-\deg\gamma$. In particular
$\sum_{t=0}^{p-1}c^t=\rho^{p-1}$.
\item[(c)] For $\alpha=\sum_{x\in C}c_xx\in\F[C]$ put $\alpha^{(-1)}=\sum_{x\in C}c_xx^{-1}$. The map
$\alpha\mapsto\alpha^{(-1)}$ is a ring automorphism of $\F[C]$, and $\ord\alpha^{(-1)}=\ord\alpha$ for
every nonzero $\alpha\in\F[C]$.
\end{itemize}
\end{lemma}

Lemma~\ref{lem:order} is well known; part (b) is the case $m=1$ of the description of the powers
of the radical of the group algebra over $\F_p$ of an elementary abelian group of order $p^m$ as
generalized Reed--Muller codes (see \cite{Ber}, \cite{Cha}, \cite{LM}). We include the short proof,
which works over any field of characteristic $p$.

\begin{proof}
(a) The map $c\mapsto Y$ gives $\F[C]\cong\F[Y]/(Y^{p^n}-1)$, and in characteristic $p$ we have
$Y^{p^n}-1=(Y-1)^{p^n}$, so $\F[C]\cong\F[Y]/((Y-1)^{p^n})=\F[\rho]/(\rho^{p^n})$. If $\alpha=\sum_{k<p^n}a_k\rho^k$ and $o$ is the least $k$ with
$a_k\ne0$, then $\alpha=\rho^ou$ where $u=\sum_{k\ge o}a_k\rho^{k-o}$ has nonzero constant term
and is therefore a unit ($\rho$ is nilpotent). If $\alpha=\rho^ou$ and $\beta=\rho^{o'}u'$, then
$\alpha\beta=\rho^{o+o'}uu'$, which is $0$ if $o+o'\ge p^n$ and has order $o+o'$ otherwise.

(b) For $0\le l\le p-1$ and $t\in\{0,\dots,p-1\}$ we have $\binom tl\equiv B_l(t)\pmod p$ with
$B_l(X)=X(X-1)\cdots(X-l+1)/l!\in\F_p[X]$, a polynomial of degree $l$ ($l!$ is invertible modulo
$p$). For $1\le R\le p-1$ the polynomials $B_0,\dots,B_{R-1}$ span the same space as
$1,X,\dots,X^{R-1}$, so $\sigma_l=0$ for all $l<R$ if and only if $\sum_t\gamma(t)t^l=0$ for all
$l<R$, if and only if $\deg\gamma\le p-1-R$ (Lemma~\ref{lem:powersum}(b)). Let
$\delta=\deg\gamma\in\{0,\dots,p-1\}$ and $M=\min\{l:\sigma_l\ne0\}$ (if $\sigma_l=0$ for
all $l\le p-2$, then $\gamma$ is a nonzero constant $\gamma_0$ and $\sigma_{p-1}=\gamma_0\ne0$,
so $M$ exists). If $\delta\le p-2$, taking $R=p-1-\delta$ gives $M\ge p-1-\delta$; if moreover
$\delta\ge1$, taking $R=p-\delta$ shows $M<p-\delta$. If $\delta=0$, then $M\ge p-1$, hence
$M=p-1$. If $\delta=p-1$, taking $R=1$ shows $M=0$. In all cases $M=p-1-\delta$. Finally
$\sum_t\gamma(t)c^t=\sum_t\gamma(t)(1+\rho)^t=\sum_{l=0}^{p-1}\sigma_l\rho^l$ is a polynomial
in $\rho$ of degree at most $p-1$, whose order is $M$. For $\gamma\equiv1$ this gives order $p-1$,
and indeed $\sum_tY^t=(Y^p-1)/(Y-1)=(Y-1)^{p-1}$ in $\F[Y]$.

(c) The map is induced by the automorphism $x\mapsto x^{-1}$ of $C$, and
$\rho^{(-1)}=c^{-1}-1=-c^{-1}\rho$. If $\alpha=\rho^ku$ with $k=\ord\alpha$ and a unit $u$, then
$\alpha^{(-1)}=\rho^k(-c^{-1})^ku^{(-1)}$, and $(-c^{-1})^ku^{(-1)}$ is a unit.
\end{proof}

\subsection{An elementary inequality}

The following inequality is used in Sections~\ref{subsec:case1} and~\ref{subsec:case2}.

\begin{lemma}\label{corl}
Let $p$ be an odd prime, let $2\le r\le p$, and let $x_1\le x_2\le\dots\le x_r$ be positive integers
with $\sum_{i=1}^rx_i^2\ge p(p+1-r)$. Put $A=x_1$ and $S=x_2+\dots+x_r$.
\begin{itemize}
\item[(a)] If $A=1$, then $S\ge p-1$, with equality only if $r=p$ and $x_1=\dots=x_r=1$.
\item[(b)] If $A\ge2$ and $r\ge3$, then $S\ge p$.
\end{itemize}
\end{lemma}

\begin{proof}
For $2\le i\le r$ we have $px_i=x_i^2+(p-x_i)x_i$, hence
\begin{equation}\label{eq:star}
pS=\sum_{i=1}^rx_i^2-A^2+\sum_{i=2}^r(p-x_i)x_i\ \ge\ p(p+1-r)-A^2+\sum_{i=2}^r(p-x_i)x_i .
\end{equation}
(a) Let $A=1$. If $x_r\ge p$, then $S\ge p$. If $x_r=p-1$ and $r\ge3$, then $S\ge p-1+(r-2)\ge p$;
if $x_r=p-1$ and $r=2$, then $\sum_ix_i^2=1+(p-1)^2<p(p-1)$, which is impossible. So we may assume
$x_i\le p-2$ for all $i$. Then $(p-x_i)x_i\ge p-1$ for $2\le i\le r$, with equality only if
$x_i=1$, and \eqref{eq:star} gives $pS\ge p(p+1-r)-1+(r-1)(p-1)=p^2-r\ge p(p-1)$. Hence $S\ge p-1$,
and $S=p-1$ forces $r=p$ and equality throughout, that is, $x_i=1$ for all $i$.

(b) Let $A\ge2$ and $r\ge3$. If $2A\ge p$, then $S\ge(r-1)A\ge2A\ge p$. If $x_r\ge p-A+1$, then
$S\ge x_r+(r-2)A\ge p+1$. So we may assume $2\le A\le(p-1)/2$, in particular $p\ge5$, and
$A\le x_i\le p-A$ for all $i$. Then $(p-x_i)x_i\ge(p-A)A$ for $2\le i\le r$, and \eqref{eq:star}
gives
\[
pS-p(p-1)\ \ge\ p(2-r)-A^2+(r-1)(p-A)A=(r-2)\bigl((p-A)A-p\bigr)+(p-2A)A .
\]
Since $(p-A)A\ge2(p-2)>p$ and $(p-2A)A>0$, we get $pS>p(p-1)$, so $S\ge p$.
\end{proof}

\section{The local-ring theorem}\label{sec:local}

Throughout this section $p$ is an arbitrary prime and $h=p^ch'$ with $c\ge0$ and $p\nmid h'$. Since
$\mu_h=\mu_{p^c}\times\mu_{h'}$, every $a\in\mu_h$ can be written uniquely as
$a=\zeta_{p^c}^{e}\eta$ with $e\in\{0,\dots,p^c-1\}$ and $\eta\in\mu_{h'}$.

\begin{theorem}\label{thm:etaconst}
Let $D=\sum_{x\in G}a_xx$ be a $\BH(G,p^ch')$ matrix with $c\ge0$ and $p\nmid h'$, and write
$a_x=\zeta_{p^c}^{a(x)}\eta(x)$ with $a(x)\in\{0,\dots,p^c-1\}$ and $\eta(x)\in\mu_{h'}$. Then
$\eta$ is constant on every coset of $H$: for every $t\in\{0,\dots,p-1\}$ the function
$y\mapsto\eta(g^{t+py})$ is constant.
\end{theorem}

\begin{proof}
Let $O=\Z[\zeta_{p^ch'}]$, fix a prime ideal $\PP$ of $O$ containing $p$, let $\F=O/\PP$ (a finite
field of characteristic $p$), and let $f\colon O\to\F$ be the reduction map. We also write $f$ for
the ring homomorphism $O[G]\to\F[G]$ given by $\sum_xc_xx\mapsto\sum_xf(c_x)x$.

\emph{(i)} First, if $c\ge1$, then $1-\zeta_{p^c}\in\PP$: indeed
$\prod_{1\le k<p^c,\ p\nmid k}(1-\zeta_{p^c}^k)=\Phi_{p^c}(1)=p\in\PP$, so some factor
$1-\zeta_{p^c}^k$ lies in $\PP$, and $1-\zeta_{p^c}$ is a multiple of it (write
$\zeta_{p^c}=(\zeta_{p^c}^k)^{k'}$ with $kk'\equiv1\pmod{p^c}$). Hence, for every $c\ge0$, we have
$f(a_x)=f(\eta(x))$ and $f(\ov{a_x})=f\bigl(\zeta_{p^c}^{-a(x)}\eta(x)^{-1}\bigr)
=f(\eta(x))^{-1}$. Second, the reduction is injective on $\mu_{h'}$: since
$\prod_{k=1}^{h'-1}(1-\zeta_{h'}^k)=h'\notin\PP$, we have
$\zeta_{h'}^u-\zeta_{h'}^v=\zeta_{h'}^u(1-\zeta_{h'}^{v-u})\notin\PP$ for $u\not\equiv v\pmod{h'}$.

\emph{(ii)} Applying $f$ to $DD^{(-1)}=p^2$, we obtain
\[
\begin{gathered}
f(D)\,f(D^{(-1)})=0\quad\text{in }\F[G],\\
f(D)=\sum_{x\in G}f(\eta(x))\,x,\qquad
f(D^{(-1)})=\sum_{x\in G}f(\eta(x))^{-1}x^{-1}.
\end{gathered}
\]
Both factors are nonzero, since their coefficients lie in $\F^*$.

\emph{(iii)} By Lemma~\ref{lem:order}(a), $\F[G]=\F[\rho]/(\rho^{p^2})$
with $\rho=g-1$. For $t\in\{0,\dots,p-1\}$ let $\gamma_t\colon\F_p\to\F^*$ with
$\gamma_t(y)=f(\eta(g^{t+py}))$ be the restriction of $f\circ\eta$ to the coset $g^tH$, and let
$d_t=\deg\gamma_t$ and $d'_t=\deg(\gamma_t^{-1})$, where $\gamma_t^{-1}$ is the pointwise inverse.
In characteristic $p$ we have $g^{t+py}=(1+\rho)^t(1+\rho)^{py}=(1+\rho)^t(1+\rho^p)^y$, and
$(1+\rho^p)^y=\sum_{l=0}^{p-1}\binom yl\rho^{pl}$ for $0\le y\le p-1$. Hence, as polynomials in
$\rho$,
\[
f(D)=\sum_{t=0}^{p-1}(1+\rho)^t\sum_{y=0}^{p-1}\gamma_t(y)(1+\rho^p)^y
=\sum_{l=0}^{p-1}\rho^{pl}\,T_l(\rho),
\]
where
\[
T_l(\rho)=\sum_{t=0}^{p-1}\sigma_{t,l}(1+\rho)^t,\qquad \sigma_{t,l}=\sum_{y}\gamma_t(y)\binom yl .
\]
Each $T_l$ has degree at most $p-1$, so the right side has degree at most $p^2-1$ and is the
reduced representative of $f(D)$ in $\F[\rho]/(\rho^{p^2})$; the coefficient of
$\rho^{pl+k}$ ($0\le l,k\le p-1$) in it is the coefficient of $\rho^k$ in $T_l$. Since
$(1+\rho)^t$ has degree exactly $t$, these polynomials are linearly independent, and $T_l\ne0$ if
and only if $\sigma_{t,l}\ne0$ for some $t$. By Lemma~\ref{lem:order}(b),
$\min\{l:\sigma_{t,l}\ne0\}=p-1-d_t$. Hence $L_0:=\min\{l:T_l\ne0\}=p-1-d$ with $d=\max_td_t$,
and
\[
\ord f(D)=pL_0+\ord_\rho T_{L_0}\le p(p-1-d)+p-1 .
\]

\emph{(iv)} Let $\beta=\sum_xf(\eta(x))^{-1}x$, whose coefficient function on the coset $g^tH$ is
$\gamma_t^{-1}$. Then $f(D^{(-1)})=\beta^{(-1)}$ in the notation of Lemma~\ref{lem:order}(c), so
$\ord f(D^{(-1)})=\ord\beta$, and the computation of (iii), applied to $\beta$, gives
$\ord f(D^{(-1)})\le p(p-1-d')+p-1$ with $d'=\max_td'_t$.

\emph{(v)} By (ii) and Lemma~\ref{lem:order}(a), $\ord f(D)+\ord
f(D^{(-1)})\ge p^2$. Hence $p^2\le p(2p-2-d-d')+2p-2$, i.e.\ $p(d+d')\le p^2-2<p^2$, so $d+d'\le p-1$.
Consequently $d_t+d'_t\le d+d'\le p-1$ for every $t$, and $\gamma_t\cdot\gamma_t^{-1}\equiv1$; by
Lemma~\ref{lem:powersum}(c), $\gamma_t$ is constant. By the injectivity in (i), $\eta$ is constant
on the coset $g^tH$.
\end{proof}

The following necessary condition was previously established in \cite[Theorem~3.6]{DS}.

\begin{corollary}\label{cor:pdividesh}
If $\BH(G,h)\ne\emptyset$, then $p\mid h$.
\end{corollary}

\begin{proof}
Suppose that $p\nmid h$ and let $D\in\BH(G,h)$. Applying Theorem~\ref{thm:etaconst}
with $c=0$ shows that the coefficients of $D$ are constant on each coset of $H$. Thus
\[
D=\left(\sum_{t=0}^{p-1}\omega_tg^t\right)\left(\sum_{y=0}^{p-1}g^{py}\right)
\]
for some $\omega_t\in\mu_h$. A character $\chi$ of $G$ that is nontrivial on $H$
annihilates the second factor, so $\chi(D)=0$, contradicting $|\chi(D)|^2=p^2$.
\end{proof}

\section{\texorpdfstring{The case $\nu_p(h)=2$ for odd $p$}{The case nu\_p(h)=2 for odd p}}\label{sec:nu2}
Recall that $H=\langle w \rangle$ denotes the subgroup of $G$ of order $p$.
For an arbitrary prime $p$ and $i,j\in\{0,\dots,p-1\}$ let $\tau_{i,j}$ be the character of $G$
with $\tau_{i,j}(g)=\zeta_{p^2}^{-(i+pj)}$, and put $\tau_i=\tau_{i,0}$. These are all the
characters of $G$. Since $\tau_{i,j}(w)=\zeta_p^{-i}$ does not depend on $j$, the characters of $H$ are
the restrictions of $\tau_0,\dots,\tau_{p-1}$ to $H$, and $\tau_{i,j}$ restricts to $\tau_i$. This
notation will also be used in Sections~\ref{sec:nu3}, \ref{sec:p2}, and~\ref{sec:proofmain}.

Throughout this section $p$ is an odd prime, $h=p^2h'$ with $p\nmid h'$, and $D=\sum_xa_xx$ is a
$\BH(G,h)$ matrix. We write $a_x=\zeta_{p^2}^{a(x)}\eta(x)$ with $a(x)\in\{0,\dots,p^2-1\}$ and
$\eta(x)\in\mu_{h'}$, and $x=g^{t(x)+py(x)}$ with $t(x),y(x)\in\{0,\dots,p-1\}$. By
Theorem~\ref{thm:etaconst} there are $\omega_t\in\mu_{h'}$ with
\begin{equation}\label{eq:etaomega}
\eta(g^{t+py})=\omega_t\qquad(0\le t,y\le p-1).
\end{equation}

\subsection{Field descent}\label{sec:descent}

Fix $i\in\{0,\dots,p-1\}$. For $x\in G$ define $s_i(x),k_i(x)\in\{0,\dots,p-1\}$ by
\begin{equation}\label{eq:sk}
s_i(x)+pk_i(x)\equiv a(x)-i\bigl(t(x)+py(x)\bigr)\pmod{p^2}.
\end{equation}
Since $\tau_{i,j}(x)=\zeta_{p^2}^{-i(t(x)+py(x))}\zeta_p^{-jt(x)}$ for $j\in\{0,\dots,p-1\}$, each
term of $\tau_{i,j}(D)$ has the form
\[
a_x\tau_{i,j}(x)=\zeta_{p^2}^{\,s_i(x)}\zeta_p^{\,k_i(x)-jt(x)}\,\eta(x),
\]
and grouping the terms of $\tau_{i,j}(D)$ according to $s_i(x)$ gives
\[
\tau_{i,j}(D)=\sum_{s=0}^{p-1}Y_{i,j,s} \text{ with }
Y_{i,j,s}=\zeta_{p^2}^s\sum_{x:\;s_i(x)=s}\zeta_p^{k_i(x)-jt(x)}\eta(x).
\]

\begin{lemma}\label{cor:unique}
For all $i,j\in\Z_p$ there is exactly one $s\in\Z_p$ with $Y_{i,j,s}\ne0$, denoted by
$s_0(i,j)$, and $|Y_{i,j,s_0(i,j)}|=p$.
\end{lemma}

\begin{proof}
We have $X\ov X=p^2$ for $X=\tau_{i,j}(D)$. By \eqref{eq:coset} with $n=b=2$, we have
$X\in\zeta_{p^2}^{e}\Q(\zeta_{ph'})$ for some $e=s_0+pk_0$ with $s_0,k_0\in\{0,\dots,p-1\}$. As
$\zeta_{p^2}^{pk_0}\in\Q(\zeta_{ph'})$, we get $X\in\zeta_{p^2}^{s_0}\Q(\zeta_{ph'})$. Comparing this
with $X=\sum_sY_{i,j,s}$, we obtain from Lemma~\ref{lem:basis}(a) that $Y_{i,j,s}=0$ for
$s\ne s_0$ and $Y_{i,j,s_0}=X$, which has absolute value $p$.
\end{proof}

\subsection{Blocks and lines}\label{subsec:setting}

The \emph{class} of $x\in G$ is the number $\sigma(x)\in\{0,\dots,p-1\}$ with
$\sigma(x)\equiv a(x)\pmod p$. The
\emph{blocks} of $D$ are
\[
D_{t,s}=\sum_{x:\;t(x)=t,\ \sigma(x)=s}\zeta_{p^2}^{a(x)}\eta(x)\,w^{y(x)}\in\Z[\zeta_{p^2h'}][H]
\qquad(t,s\in\Z_p),
\]
so that $D=\sum_{t,s=0}^{p-1}g^tD_{t,s}$.
Let $N_{t,s}$ be the number of  elements $x$ with $t(x)=t$ and $\sigma(x)=s$ (the number of terms of
$D_{t,s}$). Then $\sum_sN_{t,s}=p$ for every $t$  and $\sum_{t,s}N_{t,s}=p^2$. We picture the blocks as
the points $(t,s)$ of the affine plane $\Z_p^2$, with \emph{mass} $N_{t,s}$; a point is
\emph{occupied} if $N_{t,s}>0$. The \emph{column} of $t$ is $\{(t,s):s\in\Z_p\}$; it has mass $p$.
Let
\[
A=\min\{N_{t,s}:\ N_{t,s}>0\}.
\]

\begin{lemma}\label{lem:blocknonzero}
\begin{itemize}
\item[(a)] If $1\le N_{t,s}\le p-1$, then $\tau_i(D_{t,s})\ne0$ for every $i\in\Z_p$.
\item[(b)] Either $A=1$, or $2\le A\le(p-1)/2$, or $A=p$. If $A=p$, then every column contains
exactly one occupied point, i.e., all elements of a coset $g^tH$ have the same class.
\end{itemize}
\end{lemma}

\begin{proof}
(a) By \eqref{eq:etaomega}, $\tau_i(D_{t,s})=\omega_t\zeta_{p^2}^{s}\sum_x\zeta_{p^2}^{a(x)-s}\zeta_p^{-iy(x)}$,
the sum over the $N_{t,s}$ elements $x$ of the block. Since $a(x)-s\equiv0\pmod p$, every term of the
sum is a $p$-th root of unity, and the sum is nonzero by Lemma~\ref{lem:sums}.
(b) If $A\ge(p+1)/2$, then every column, which has mass $p$, contains exactly one occupied point,
since two occupied points would have total mass at least $2A>p$; that point has mass $p$, so $A=p$.
\end{proof}

\begin{definition}
For $i,s\in\Z_p$, 
define the \emph{line} $L_{i,s}=\{(t,it+s):\ t\in\Z_p\}$,
its \emph{slope} as $i$,
and its \emph{mass} as
$m(L_{i,s})=\sum_{t\in \Z_p}N_{t,it+s}$.
\end{definition}

A line meets every column in exactly one point; two distinct lines through a point have different
slopes and meet only in that point; the $p$ lines of slope $i$ partition $\Z_p^2$. 
Note that $(t(x),\sigma(x))\in L_{i,s}$ if and only if $\sigma(x)\equiv it(x)+s \pmod p$.
Moreover, $s_i(x)\equiv a(x)-it(x)\equiv\sigma(x)-it(x)\pmod p$ by \eqref{eq:sk}.
Thus
\begin{equation}\label{eq:online}
(t(x),\sigma(x))\in L_{i,s}\quad\text{if and only if}\quad s_i(x)=s .
\end{equation}
Hence
\[
Y_{i,j,s}=\zeta_{p^2}^s\sum_{x:\;(t(x),\sigma(x))\in L_{i,s}}\zeta_p^{k_i(x)-jt(x)}\eta(x)
\]
is the contribution of the line $L_{i,s}$ to $\tau_{i,j}(D)$, and by Lemma~\ref{cor:unique}
exactly one line of slope $i$ contributes to $\tau_{i,j}(D)$. For a line $L=L_{i,s}$ let
\[
J(L)=\{j\in\Z_p:\ Y_{i,j,s}\ne0\}=\{j\in\Z_p:\ s_0(i,j)=s\}
\]
be the set of $j$ for which this line is $L$. For fixed $i$ the sets $J(L_{i,s})$ with $s\in\Z_p$ partition $\Z_p$.
We call $L$ \emph{active} if $J(L)\ne\emptyset$, \emph{inactive} otherwise, and \emph{flat} if
$J(L)=\Z_p$.

\subsection{Mass bounds for lines}\label{subsec:lines}

Fix a slope $i\in\{0,\dots,p-1\}$. For a line $L=L_{i,s}$ put
\[
W_L=\sum_{t=0}^{p-1}W_{L,t}\,w^t\in\Z[\zeta_{p^2h'}][H]\qquad\text{with}\qquad
W_{L,t}=\zeta_{p^2}^{-it}\tau_i(D_{t,it+s})=\tau_i(g^tD_{t,it+s}).
\]

\begin{lemma}\label{lem:lineid}
Let $L=L_{i,s}$ be a line.
\begin{itemize}
\item[(a)] We have $W_{L,t}=\sum_{x:\;t(x)=t,\ s_i(x)=s}\zeta_{p^2}^{s+pk_i(x)}\eta(x)$ and
$\tau_j(W_L)=Y_{i,j,s}$ for every $j\in\Z_p$.
\item[(b)] We have $W_LW_L^{(-1)}=p\,\Theta_{J(L)}$, where
$\Theta_J:=\sum_{d\in\Z_p}\bigl(\sum_{j\in J}\zeta_p^{jd}\bigr)w^d$ for $J\subseteq\Z_p$. In
particular $\sum_t|W_{L,t}|^2=p\,|J(L)|$, and $W_LW_L^{(-1)}=p^2$ if $L$ is flat.
\item[(c)] If $L$ is inactive, then $\tau_i(D_Q)=0$ for every point $Q$ of $L$.
\end{itemize}
\end{lemma}

\begin{proof}
(a) By \eqref{eq:online}, the $x$ with $t(x)=t$ and $s_i(x)=s$ are exactly the elements of the block
$D_{t,it+s}$, and for them
\[
\zeta_{p^2}^{-it}\zeta_{p^2}^{a(x)}\eta(x)\zeta_p^{-iy(x)}=\zeta_{p^2}^{a(x)-i(t+py(x))}\eta(x)
=\zeta_{p^2}^{s+pk_i(x)}\eta(x)
\]
by \eqref{eq:sk}. Summing over the block gives the first claim. Hence
\[
\tau_j(W_L)=\sum_tW_{L,t}\zeta_p^{-jt}=\sum_{x:\,s_i(x)=s}\zeta_{p^2}^{s}\zeta_p^{k_i(x)-jt(x)}\eta(x)=Y_{i,j,s}.
\]
(b) By (a) and Lemma~\ref{cor:unique}, $\tau_{j'}(W_LW_L^{(-1)})=|Y_{i,j',s}|^2$ equals $p^2$ if
$j'\in J(L)$ and $0$ otherwise, and $\tau_{j'}(p\Theta_{J(L)})=p\sum_{j\in J(L)}\sum_d\zeta_p^{(j-j')d}$
takes the same values. By Fourier inversion on $H$ the two elements agree. Comparing the coefficients of $w^0$ gives
$\sum_t|W_{L,t}|^2=p|J(L)|$, and $\Theta_{\Z_p}=p$.
(c) If $J(L)=\emptyset$, then all $W_{L,t}$ vanish by (b), and (a) gives the claim.
\end{proof}

\begin{lemma}\label{lem:line}
Let $P=(t_1,s_1)$ be a point with $0<N_P<p$, let $L=L_{i,s}$ be a line through $P$, and for
$t\in\Z_p$ let $x_t$ be the mass of the point $(t,it+s)$ of $L$. Let
$T=\{t\in\Z_p:\ W_{L,t}\ne0\}$. Then $t_1\in T$ and $J(L)\ne\emptyset$; moreover $x_t\ge|W_{L,t}|$ for
all $t$, and
\[
m(L)\ \ge\ \sum_{t\in T}x_t,\qquad \sum_{t\in T}x_t^2\ \ge\ p\,|J(L)| .
\]
Moreover:
\begin{itemize}
\item[(a)] If $L$ is flat, then $|T|\ge3$, and $|T|\ge4$ if $p\ge5$.
\item[(b)] If $L$ is not flat, then $|T|\ge p+1-|J(L)|$ and $|T|(|T|-1)+1\ge p$.
\end{itemize}
\end{lemma}

\begin{proof}
By the definition of $W_L$, $|W_{L,t}|=|\tau_i(D_{t,it+s})|$, and $\tau_i(D_{t,it+s})$ is a sum of
$x_t$ roots of unity, so $x_t\ge|W_{L,t}|$. As $0<N_P<p$, we have $\tau_i(D_P)\ne0$ by
Lemma~\ref{lem:blocknonzero}(a), so $t_1\in T$. By Lemma~\ref{lem:lineid}(b),
$\sum_{t\in T}x_t^2\ge\sum_t|W_{L,t}|^2=p|J(L)|$, and $J(L)\ne\emptyset$ as $W_{L,t_1}\ne0$. Clearly
$m(L)\ge\sum_{t\in T}x_t$.

(a) Let $L$ be flat, so that $W_LW_L^{(-1)}=p^2$ by Lemma~\ref{lem:lineid}(b). For $d\ne0$ the
coefficient of $w^d$ in $W_LW_L^{(-1)}$ is $\sum_tW_{L,t+d}\ov{W_{L,t}}$, a sum over the pairs
$(t+d,t)\in T\times T$, and it vanishes. Hence no $d\ne0$ is the difference $t-t'$ of exactly
one pair $(t,t')\in T\times T$, since otherwise the sum for $d$ would consist of a single nonzero
product; we say that $T$ has no \emph{unique differences}. The case $|T|=1$ is impossible, as then
$p^2=\sum_t|W_{L,t}|^2=|W_{L,t_1}|^2\le N_P^2<p^2$. The case $|T|=2$ is impossible, since for
$T=\{a,b\}$ the difference $b-a$ is unique ($b-a=a-b$ would force $2(b-a)=0$, but $p$ is odd).
Suppose $T=\{a,b,c\}$ and $p\ge5$. The difference $b-a$ must coincide with one of
$c-a,\ a-b,\ b-c,\ c-b,\ a-c$; the first three are impossible (the first forces $b=c$, the second
$2(b-a)=0$, the third $a=c$), so $b-a\in\{c-b,\,a-c\}$, and in both cases $T$ is an arithmetic
progression $\{m_0-d,m_0,m_0+d\}$ with $d\ne0$. Its difference $2d=(m_0+d)-(m_0-d)$ occurs for
another pair only if $2d\in\{d,-d,-2d\}$; as $d\ne0$ and $p$ is odd, this forces $3d=0$, i.e.\
$p=3$, a contradiction. Thus $|T|\ge4$ if $p\ge5$.

(b) Here $1\le|J(L)|\le p-1$, so by Lemma~\ref{lem:lineid}(b) and Lemma~\ref{lem:sums} the
coefficient $p\sum_{j\in J(L)}\zeta_p^{jd}$ of $w^d$ in $W_LW_L^{(-1)}$ is nonzero for every $d\ne0$.
As in (a), this coefficient is a sum over the pairs $(t+d,t)\in T\times T$, so every $d\ne0$ is a
difference of two elements of $T$, and $|T|(|T|-1)+1\ge p$. By Lemma~\ref{lem:lineid}(a), $\tau_j(W_L)\ne0$ if and only if
$j\in J(L)$, so Lemma~\ref{lem:uncertainty} gives $|T|+|J(L)|\ge p+1$.
\end{proof}

\subsection{\texorpdfstring{The case $A=1$}{The case A=1}}\label{subsec:case1}

\begin{theorem}\label{thm:case1}
Let $p$ be odd and $D$ a $\BH(G,p^2h')$ matrix. Then no block of $D$ has exactly one term, i.e.\
$A\ne1$.
\end{theorem}

\begin{proof}
\emph{Step 1.} Let $P=(t_1,s_1)$ with $N_P=1$. For $i\in\Z_p$ let $L^{(i)}$ be the line of slope $i$
through $P$, let the set $T$ and the masses $x_t$ be as in Lemma~\ref{lem:line}, and put $n=|J(L^{(i)})|$.
Then $x_t\ge1$ for $t\in T$, $x_{t_1}=1$ and $\sum_{t\in T}x_t^2\ge pn$. If $L^{(i)}$ is flat, then
$n=p$ and $|T|\ge3$ by Lemma~\ref{lem:line}(a); otherwise $n\ge p+1-|T|$ by Lemma~\ref{lem:line}(b). In both cases $|T|\ge2$ and
$\sum_{t\in T}x_t^2\ge p(p+1-|T|)$, so Lemma~\ref{corl}(a) with $r=|T|$ gives
$m(L^{(i)})\ge\sum_{t\in T}x_t\ge p$, with equality only if $|T|=p$ and $x_t=1$ for all $t\in T$;
in that case $pn\le\sum_{t\in T}x_t^2=p$, so $n=1$.
The column of $P$ and the $p$ lines $L^{(i)}$ cover $\Z_p^2$ and
pairwise meet only in $P$, so
\[
p^2=p+\sum_{i\in\Z_p}\bigl(m(L^{(i)})-1\bigr)\ \ge\ p+p(p-1)=p^2 .
\]
Hence $m(L^{(i)})=p$ for every $i$; consequently $L^{(i)}$ is not flat, $|J(L^{(i)})|=1$ and
$|T|=p$.

\emph{Step 2.} All points of these lines are blocks with one term. Fix $i$, let
$L=L^{(i)}=L_{i,s}$, and let $J(L)=\{j\}$. As $|T|=p$, all $p$ points of $L$ are occupied, and as
$m(L)=p$, each of them is a block with exactly one term; let $z_t$ be the element of the point of $L$
in the column $t$. By Lemma~\ref{lem:line}, $|W_{L,t}|\le1$ for all $t$, and
$\sum_t|W_{L,t}|^2=p$ by Lemma~\ref{lem:lineid}(b), so $|W_{L,t}|=1$. As $\tau_{j'}(W_L)=0$ for
$j'\ne j$ (Lemma~\ref{lem:lineid}(a)), Fourier inversion on $H$ gives
$W_{L,t}=p^{-1}\tau_j(W_L)\zeta_p^{jt}$. Hence
$\zeta_p^{-jt}W_{L,t}=\zeta_{p^2}^s\zeta_p^{k_i(z_t)-jt}\eta(z_t)$ does not depend on $t$. As
$\mu_p\cap\mu_{h'}=\{1\}$, all $x$ with $s_i(x)=s$ have the same value of $k_i(x)-jt(x)$ modulo
$p$.

\emph{Step 3.} All blocks have one term. By Step~2 every point outside the column $t_1$ lies on a
line through $P$ and is therefore a block with one term. Applying Steps~1--2 to such a block $P'$
in a column $t'\ne t_1$ shows that the points of the column $t_1$ are blocks with one term as well.
Hence $N_{t,s}=1$ for all $(t,s)$: every coset $g^tH$ contains exactly one element of each class.

\emph{Step 4.} Now every line has mass $p$, and every point is a block with one
term, so by Step~2 (applied to any point of the line) every line $L$ satisfies $|J(L)|=1$. For
fixed $i$ the $p$ lines of slope $i$ therefore correspond bijectively to the elements $j\in\Z_p$;
write $j=j(i,s)$ for the line $L_{i,s}$, and $\kappa(i,s)\in\Z_p$ for the common value of
$k_i(x)-j(i,s)t(x)$ on it (Step~2). Thus, for all $x=g^{t}w^{y}$ with $s_i(x)=s$,
\begin{equation}\label{eq:lin}
a(x)-(i+pj(i,s))(t+py)\equiv s+p\kappa(i,s)\pmod{p^2},
\end{equation}
and $s\mapsto j(i,s)$ is a permutation of $\Z_p$ for every $i$.

\emph{Step 5.} Let $x=g^tw^y$ and $v=\sigma(x)$. For $i=0$ the line of $x$ is
$L_{0,v}$, and \eqref{eq:lin} reads $a(x)\equiv v+p\bigl(\kappa(0,v)+j(0,v)t\bigr)\pmod{p^2}$. For
$i=1$ the line of $x$ is $L_{1,u}$ with $u\equiv v-t$; write $v-t=u-p\delta(v,t)$ with
$u\in\{0,\dots,p-1\}$ and $\delta(v,t)\in\{0,1\}$. Substituting $a(x)$ into \eqref{eq:lin} for
$i=1$ and dividing by $p$, we get
\[
y\equiv\kappa(0,v)-\kappa(1,u)-\delta(v,t)+\bigl(j(0,v)-j(1,u)\bigr)t\pmod p .
\]
Fix $t$ and let $v$ run through $\{0,\dots,p-1\}$. By Step~3, $y=y(t,v)$ then runs through
$\{0,\dots,p-1\}$, and so does $u$. Summing over $v$, and using that $s\mapsto j(0,s)$ and
$s\mapsto j(1,s)$ are permutations and that $p$ is odd (so that $\sum_{r=0}^{p-1}r\equiv0\pmod p$),
we get $0\equiv K-\sum_{v}\delta(v,t)=K-t\pmod p$, where $K=\sum_s\kappa(0,s)-\sum_s\kappa(1,s)$
does not depend on $t$, and where we used that $\delta(v,t)=1$ if and only if $v<t$. Taking
$t=0$ and $t=1$ gives $0\equiv1\pmod p$, a contradiction.
\end{proof}

\subsection{\texorpdfstring{The case $2\le A\le(p-1)/2$}{The case 2<=A<=(p-1)/2}}\label{subsec:case2}

Here $p\ge5$, since for $p=3$ this range is empty.

\begin{theorem}\label{thm:case2}
Let $p\ge5$ and let $D$ be a $\BH(G,p^2h')$ matrix. Then $A\notin[2,(p-1)/2]$.
\end{theorem}

\begin{proof} Suppose that $2\le A\le(p-1)/2$, and let $P=(t_1,s_1)$ be a point with $N_P=A$. Let $L$ be a line
through $P$, and let $T$ and $x_t$ be as in Lemma~\ref{lem:line}. Then $x_t\ge A$ for $t\in T$, $x_{t_1}=A$
and $\sum_{t\in T}x_t^2\ge p|J(L)|$. If $L$ is flat, then $|J(L)|=p$ and $|T|\ge4$ by
Lemma~\ref{lem:line}(a). If $L$ is not flat, then $|J(L)|\ge p+1-|T|$, and $|T|\ge3$ as
$|T|(|T|-1)+1\ge p\ge5$ by Lemma~\ref{lem:line}(b). In both cases $|T|\ge3$ and
$\sum_{t\in T}x_t^2\ge p(p+1-|T|)$, so Lemma~\ref{corl}(b) gives
$m(L)-A\ge\sum_{t\in T}x_t-A\ge p$. As the column of $P$ and the sets $L\setminus\{P\}$ for the $p$
lines $L$ through $P$ are pairwise disjoint, we get $p^2\ge p+p\cdot p$, a contradiction.
\end{proof}

\subsection{\texorpdfstring{The case $A=p$}{The case A=p}}\label{subsec:case3}

In this subsection we assume $A=p$. By Lemma~\ref{lem:blocknonzero}(b), every coset $g^tH$
consists of elements of a single class $s_t$; so $a(g^{t+py})=s_t+p\,b_t(y)$ with
$s_t,b_t(y)\in\{0,\dots,p-1\}$, and by \eqref{eq:etaomega}
\[
a_{g^{t+py}}=\zeta_{p^2}^{s_t}\,\omega_t\,\zeta_p^{b_t(y)},\qquad
D=\sum_{t=0}^{p-1}g^t D_{t,s_t},\qquad D_{t,s_t}=\zeta_{p^2}^{s_t}\sum_{y=0}^{p-1}f_t(y)w^y,\quad
f_t(y)=\omega_t\zeta_p^{b_t(y)} .
\]
We put
\[
\fh_t(i):=\sum_yf_t(y)\zeta_p^{-iy}=\zeta_{p^2}^{-s_t}\tau_i(D_{t,s_t})=\omega_tX_{t,i},\qquad
X_{t,i}=\sum_{y}\zeta_p^{b_t(y)-iy}=\sum_{k=0}^{p-1}n^i_t(k)\zeta_p^k,
\]
where $n^i_t(k)=|\{y:\ b_t(y)-iy\equiv k\pmod p\}|$; thus $n^i_t\ge0$ and $\sum_kn^i_t(k)=p$.
Since $\tau_{i,j}(g^tX)=\zeta_{p^2}^{-(i+pj)t}\tau_i(X)$ for $X\in\C[H]$, we have
\begin{equation}\label{eq:tau3}
\tau_{i,j}(D)=\sum_{t=0}^{p-1}\zeta_{p^2}^{s_t-it}\zeta_p^{-jt}\,\fh_t(i)
=\tau_j\Bigl(\sum_{t=0}^{p-1}\zeta_{p^2}^{s_t-it}\fh_t(i)\,w^t\Bigr).
\end{equation}
Parseval's identity on $H$, applied to $\sum_yf_t(y)w^y$ (whose $p$ coefficients have absolute
value $1$) and to the element of $\C[H]$ in \eqref{eq:tau3} (whose character values $\tau_{i,j}(D)$ have
absolute value $p$), gives
\begin{equation}\label{eq:parseval3}
\sum_{i\in\Z_p}|\fh_t(i)|^2=p^2\quad\text{for every }t,\qquad
\sum_{t\in\Z_p}|\fh_t(i)|^2=p^2\quad\text{for every }i .
\end{equation}

\paragraph{Lines in the case $A=p$.} The only occupied point of the column $t$ is $(t,s_t)$. For a
line $L$ let $T(L)=\{t:(t,s_t)\in L\}$. For fixed $i$ the sets $T(L_{i,s})$ with $s\in\Z_p$ partition
$\Z_p$, and two lines of different slopes have at most one point in common. For $L=L_{i,s}$ we have
$W_{L,t}=\zeta_{p^2}^{-it}\tau_i(D_{t,s_t})=\zeta_{p^2}^{s_t-it}\fh_t(i)$ if $t\in T(L)$, and
$W_{L,t}=0$ otherwise, as the point of $L$ in the column $t$ is then unoccupied. Hence
Lemma~\ref{lem:lineid}(b) gives
\begin{equation}\label{eq:lineid3}
\sum_{t\in T(L)}|\fh_t(i)|^2=p\,|J(L)| .
\end{equation}

\begin{proposition}\label{prop:oneline}
For every $i\in\Z_p$ exactly one line $L_i$ of slope $i$ is active, and it is flat. Moreover
$\fh_t(i)=0$ for all $t\notin T(L_i)$.
\end{proposition}

\begin{proof}
For an active line $L=L_{i,s}$ the set of $t$ with $W_{L,t}\ne0$ is
$T'(L)=\{t\in T(L):\ \fh_t(i)\ne0\}$, and by Lemma~\ref{lem:lineid}(a) we have $\tau_j(W_L)\ne0$
if and only if $j\in J(L)$. So Lemma~\ref{lem:uncertainty} gives $|T'(L)|+|J(L)|\ge p+1$. Let $m$
be the number of active lines of slope $i$. As the sets $J(L_{i,s})$ with $s\in\Z_p$ partition
$\Z_p$, we have $m\ge1$, and summing over the active lines of slope $i$ gives
$\sum_L|T'(L)|\ge m(p+1)-p$. The sets $T'(L)$ are pairwise disjoint subsets of $\Z_p$, so
$m(p+1)\le2p$ and $m=1$. Hence exactly one line $L_i$ of slope $i$ is active, and $J(L_i)=\Z_p$. If
$t\notin T(L_i)$, then $(t,s_t)$ lies on an inactive line of slope $i$, and $\fh_t(i)=0$ by
\eqref{eq:lineid3}.
\end{proof}

For $t,i\in\Z_p$ let $\ell_t(i)=\max_kn^i_t(k)$. By Lemma~\ref{lem:sums}, $\fh_t(i)=0$ if and only
if $n^i_t\equiv1$, that is, if and only if $y\mapsto b_t(y)-iy$ is a permutation of $\Z_p$; and
$|\fh_t(i)|=p$ if and only if $\ell_t(i)=p$, as $X_{t,i}$ is a sum of $p$ roots of unity. The
number $\sum_kn^i_t(k)^2-p$ counts the ordered pairs $(y,y')$ with $y\ne y'$ and
$b_t(y)-b_t(y')\equiv i(y-y')\pmod p$, so the next lemma says that each $i$ arises from exactly
$p(p-1)$ of the $p^2(p-1)$ such pairs with $t\in\Z_p$.

\begin{lemma}\label{lem:pairs}
For every $i\in\Z_p$ we have $\sum_{t\in\Z_p}\bigl(\sum_kn^i_t(k)^2-p\bigr)=p(p-1)$.
\end{lemma}

\begin{proof}
For $c\in\Z_p$ put $X^{(c)}_{t,i}=\sum_kn^i_t(k)\zeta_p^{ck}$, so that $X^{(0)}_{t,i}=p$ and
$X^{(1)}_{t,i}=X_{t,i}$. These are the character values of $\sum_kn^i_t(k)w^k\in\Z[H]$, so
Parseval's identity gives $\sum_c|X^{(c)}_{t,i}|^2=p\sum_kn^i_t(k)^2$. For $c\ne0$ the automorphism
$\zeta_p\mapsto\zeta_p^c$ of $\Q(\zeta_p)$ commutes with complex conjugation and maps
$\sum_t|X_{t,i}|^2=\sum_t|\fh_t(i)|^2=p^2$ (see \eqref{eq:parseval3}) to
$\sum_t|X^{(c)}_{t,i}|^2=p^2$. Summing over $t$ and $c$ gives
$p\sum_t\sum_kn^i_t(k)^2=p\cdot p^2+(p-1)p^2$, which is the claim.
\end{proof}

\begin{lemma}\label{lem:directions}
Let $t,i\in\Z_p$ with $\ell_t(i)\le p-1$. Then $\fh_t(a)\ne0$ for at least $\ell_t(i)$ elements
$a\in\Z_p\setminus\{i\}$.
\end{lemma}

\begin{proof}
Let $\ell=\ell_t(i)$, and let $y_1,\dots,y_\ell$ be distinct with $b_t(y_k)-iy_k\equiv k_0\pmod p$
for all $k$. As $\ell\le p-1$, there is $y\notin\{y_1,\dots,y_\ell\}$. For $1\le k\le\ell$ let
$a_k\in\Z_p$ be the solution of $a_k(y-y_k)\equiv b_t(y)-b_t(y_k)\pmod p$, so that
$b_t(y)-a_ky\equiv b_t(y_k)-a_ky_k$. Then $a_k\ne i$, since $b_t(y)-iy\not\equiv k_0$. If
$a_k=a_{k'}=a$, then $b_t(y_k)-b_t(y_{k'})$ is congruent to both $a(y_k-y_{k'})$ and
$i(y_k-y_{k'})$, which forces $y_k=y_{k'}$ as $a\ne i$. Hence $a_1,\dots,a_\ell$ are distinct. For
every $k$ the map $y'\mapsto b_t(y')-a_ky'$ is not injective, so $\fh_t(a_k)\ne0$.
\end{proof}

\begin{proposition}\label{prop:claim}
For every $i\in\Z_p$ there is exactly one $t\in\Z_p$ with $\fh_t(i)\ne0$, and $|\fh_t(i)|=p$ for
this $t$.
\end{proposition}

\begin{proof}
Fix $i$ and let $T_i=\{t:\ \fh_t(i)\ne0\}$, which is nonempty by \eqref{eq:parseval3}. Suppose that
$|T_i|\ge2$. If $\ell_t(i)=p$ for some $t\in T_i$, then $|\fh_t(i)|^2=p^2$, and \eqref{eq:parseval3}
gives $\fh_{t'}(i)=0$ for all $t'\ne t$, which is impossible. So $\ell_t(i)\le p-1$ for all
$t\in T_i$. For $t\notin T_i$ we have $n^i_t\equiv1$, and for $t\in T_i$ we have
$\sum_kn^i_t(k)^2\le\ell_t(i)\sum_kn^i_t(k)=p\,\ell_t(i)$. Therefore Lemma~\ref{lem:pairs} gives
$\sum_{t\in T_i}(\ell_t(i)-1)\ge p-1$, that is, $\sum_{t\in T_i}\ell_t(i)\ge p-1+|T_i|>p-1$. By
Lemma~\ref{lem:directions}, the sets $\{a\ne i:\ \fh_t(a)\ne0\}$ with $t\in T_i$ are subsets of
$\Z_p\setminus\{i\}$ of total size larger than $p-1$. So there are $a\ne i$ and $t\ne t'$ in $T_i$
with $\fh_t(a)\ne0$ and $\fh_{t'}(a)\ne0$. By Proposition~\ref{prop:oneline}, both $t$ and $t'$
lie in $T(L_i)\cap T(L_a)$, so the two distinct points $(t,s_t)$ and $(t',s_{t'})$ lie on both
$L_i$ and $L_a$, two lines of different slopes. This contradiction shows $|T_i|=1$, and then
$|\fh_t(i)|=p$ for the element $t$ of $T_i$ by \eqref{eq:parseval3}.
\end{proof}

\begin{theorem}\label{thm:case3}
Let $p$ be odd and let $D$ be a $\BH(G,p^2h')$ matrix with $A=p$. Then $D$ is of the form
\eqref{eq:family}.
\end{theorem}

\begin{proof}
By Proposition~\ref{prop:claim}, $|\fh_t(i)|\in\{0,p\}$ for all $t,i$, and $|\fh_t(i)|=p$ holds if
and only if all $p$ terms $\zeta_p^{b_t(y)-iy}$ of $X_{t,i}$ are equal, that is, if and only if
$b_t(y)-iy$ is constant modulo $p$ as a function of $y$. By \eqref{eq:parseval3}, for every $t$
there is exactly one $i$ with $|\fh_t(i)|=p$; we write it as $i=\pi(t)$, so that
$b_t(y)\equiv\beta_t+\pi(t)y\pmod p$ for some constant $\beta_t$. By Proposition~\ref{prop:claim},
for every $i$ exactly one $t$ has $\pi(t)=i$: the map $\pi$ is a permutation of $\Z_p$. Finally
\[
g^tD_{t,s_t}=\zeta_{p^2}^{s_t+p\beta_t}\omega_t\,g^t\sum_y\zeta_p^{\pi(t)y}w^y,
\]
which is \eqref{eq:family} with $\omega'_t=\zeta_{p^2}^{s_t+p\beta_t}\omega_t\in\mu_h$ and the
permutation $\pi$.
\end{proof}

\subsection{\texorpdfstring{Conclusion for odd $p$}{Conclusion for odd p}}

\begin{theorem}\label{thm:nu2odd}
Let $p$ be an odd prime and $h=p^2h'$ with $p\nmid h'$. Then every $\BH(G,h)$ matrix is of the
form \eqref{eq:family}.
\end{theorem}

\begin{proof}
By Lemma~\ref{lem:blocknonzero}(b), $A\in\{1\}\cup[2,(p-1)/2]\cup\{p\}$. The first two
possibilities are excluded by Theorems~\ref{thm:case1} and~\ref{thm:case2}, so $A=p$, and
Theorem~\ref{thm:case3} applies.
\end{proof}

\section{\texorpdfstring{The case $\nu_p(h)\ge3$ for odd $p$}{The case nu\_p(h)>=3 for odd p}}\label{sec:nu3}

We reduce this case to $\nu_p(h)=2$.

\begin{theorem}\label{thm:nu3}
Let $p$ be an odd prime and $h=p^ch'$ with $c\ge3$ and $p\nmid h'$. Then every $\BH(G,h)$ matrix
is of the form \eqref{eq:family}.
\end{theorem}

\begin{proof}
Let $D=\sum_xa_xx\in\BH(G,h)$. For $x\in G$ write $a_x=\zeta_{p^c}^{e(x)}\theta(x)$ with
$0\le e(x)<p^c$ and $\theta(x)\in\mu_{h'}$, and $e(x)=\alpha(x)+p^{c-2}k(x)$ with
$0\le\alpha(x)<p^{c-2}$ and $0\le k(x)<p^2$. Then
\[
a_x=\zeta_{p^c}^{\alpha(x)}b_x,\qquad b_x=\zeta_{p^2}^{k(x)}\theta(x)\in\mu_{p^2h'} .
\]
Let $\alpha_1<\dots<\alpha_r$ be the distinct values of $\alpha$, and
$D_m=\sum_{x:\,\alpha(x)=\alpha_m}b_xx\in\Z[\zeta_{p^2h'}][G]$. Then
$D=\sum_{m=1}^r\zeta_{p^c}^{\alpha_m}D_m$, and the supports of the $D_m$ partition $G$.

\emph{Step 1.} For every character $\tau$ of $G$ there is exactly one $m$ with $\tau(D_m)\ne0$,
and $|\tau(D_m)|=p$ for it. Let $X=\tau(D)\in\Z[\zeta_{p^ch'}]$; then $X\ov X=p^2$. By
\eqref{eq:coset}, applied with $n=c$ and $b=2$, we have $X\in\zeta_{p^c}^{e}\Q(\zeta_{ph'})$ for some $e$. Since
$\zeta_{p^c}^{p^{c-2}}=\zeta_{p^2}\in\Q(\zeta_{p^2h'})$, we get $X\in\zeta_{p^c}^{e_2}\Q(\zeta_{p^2h'})$
with $0\le e_2<p^{c-2}$ the residue of $e$ modulo $p^{c-2}$. On the other hand $X=\sum_m\zeta_{p^c}^{\alpha_m}\tau(D_m)$ with
$\tau(D_m)\in\Q(\zeta_{p^2h'})$ and pairwise distinct $\alpha_m\in[0,p^{c-2})$. Since the elements
$\zeta_{p^c}^{e}$ with $0\le e<p^{c-2}$ form a basis of $\Q(\zeta_{p^ch'})$ over $\Q(\zeta_{p^2h'})$
(Lemma~\ref{lem:basis}(a)), we get $\tau(D_m)=0$ whenever $\alpha_m\ne e_2$. As $X\ne0$, there is a
(unique) $m$ with $\alpha_m=e_2$, and $\tau(D_m)=\zeta_{p^c}^{-e_2}X$ has absolute value $p$.

\emph{Step 2.} We have $D_mD_{m'}^{(-1)}=0$ for $m\ne m'$, and $\sum_mD_mD_m^{(-1)}=p^2$. For $m\ne m'$ and
every character $\tau$ we have $\tau(D_mD_{m'}^{(-1)})=\tau(D_m)\ov{\tau(D_{m'})}=0$ by Step~1, so
$D_mD_{m'}^{(-1)}=0$ by Fourier inversion. Hence
$p^2=DD^{(-1)}=\sum_{m,m'}\zeta_{p^c}^{\alpha_m-\alpha_{m'}}D_mD_{m'}^{(-1)}=\sum_mD_mD_m^{(-1)}$.

\emph{Step 3.} Choose a prime $q$ with $q>p^{c-2}$ and
$q\nmid h$ (in particular $q\ne p$), and put
\[
D'=\sum_{m=1}^r\zeta_q^{\alpha_m}D_m=\sum_{x\in G}\zeta_q^{\alpha(x)}b_x\,x .
\]
Its coefficients lie in $\mu_{p^2h'q}$, and by Step~2 we have
\[
D'D'^{(-1)}=\sum_{m,m'}\zeta_q^{\alpha_m-\alpha_{m'}}D_mD_{m'}^{(-1)}=\sum_mD_mD_m^{(-1)}=p^2 .
\]
So $D'$ is
a $\BH(G,p^2h'q)$ matrix, and $p\nmid h'q$. By Theorem~\ref{thm:nu2odd},
\[
D'=\sum_{t=0}^{p-1}\Omega_t\,g^t\sum_{y=0}^{p-1}\zeta_p^{\pi(t)y}g^{py}
\]
with $\Omega_t\in\mu_{p^2h'q}$ and a permutation $\pi$ of $\Z_p$.

\emph{Step 4.} Since $q\nmid p^2h'$, we have
$\mu_{p^2h'q}=\mu_q\times\mu_{p^2h'}$, and we can write $\Omega_t=\zeta_q^{\gamma_t}\nu_t$ with
$0\le\gamma_t<q$ and $\nu_t\in\mu_{p^2h'}$. For $x=g^{t+py}$, comparing the coefficients of $x$
in the two expressions for $D'$ gives
$\zeta_q^{\alpha(x)}b_x=\zeta_q^{\gamma_t}\cdot\nu_t\zeta_p^{\pi(t)y}$ with $b_x$ and
$\nu_t\zeta_p^{\pi(t)y}$ in $\mu_{p^2h'}$. By the uniqueness of the decomposition in
$\mu_q\times\mu_{p^2h'}$, we get $\zeta_q^{\alpha(x)}=\zeta_q^{\gamma_t}$ and $b_x=\nu_t\zeta_p^{\pi(t)y}$.
As $0\le\alpha(x)<p^{c-2}<q$ and $0\le\gamma_t<q$, the first equation gives $\alpha(x)=\gamma_t$;
in particular $\alpha(x)$ depends only on the coset $t$ of $x$. Therefore
\[
a_{g^{t+py}}=\zeta_{p^c}^{\gamma_t}\nu_t\,\zeta_p^{\pi(t)y},\qquad\text{i.e.}\qquad
D=\sum_{t=0}^{p-1}\bigl(\zeta_{p^c}^{\gamma_t}\nu_t\bigr)g^t\sum_{y=0}^{p-1}\zeta_p^{\pi(t)y}g^{py},
\]
with $\zeta_{p^c}^{\gamma_t}\nu_t\in\mu_h$. This is \eqref{eq:family}.
\end{proof}

\section{The prime 2}\label{sec:p2}

For $p=2$ we do not use field descent. Instead, the following elementary proposition, which only
uses that the entries have absolute value $1$, settles the case $p=2$ for all $h$.

\begin{proposition}\label{prop:p2all}
Let $G=\langle g\rangle$ be cyclic of order $4$, and let $D=\sum_{x\in G}a_xx\in\C[G]$ with
$|a_x|=1$ for all $x$ and $|\tau_{i,j}(D)|^2=4$ for all $i,j\in\{0,1\}$. For $t\in\{0,1\}$ put
$\rho_t=a_{g^{t+2}}/a_{g^t}$. Then $\{\rho_0,\rho_1\}=\{1,-1\}$.
\end{proposition}

\begin{proof}
Put $\omega'_t=a_{g^t}$, so that $D=\sum_{t=0}^1\omega'_tg^t(1+\rho_tg^2)$ and, as
$\tau_{i,j}(g^{t+2y})=\zeta_4^{-(i+2j)t}(-1)^{iy}$,
\[
\tau_{i,j}(D)=\omega'_0\bigl(1+(-1)^i\rho_0\bigr)+(-1)^j\zeta_4^{-i}\omega'_1\bigl(1+(-1)^i\rho_1\bigr).
\]
Choose $\xi_t\in\C$ with $\xi_t^2=\rho_t$; then $|\xi_t|=1$. Write $C_t=\mathrm{Re}(\xi_t)$ and
$S_t=\mathrm{Im}(\xi_t)$, so that $C_t^2+S_t^2=1$. Then $1+\rho_t=\xi_t(\ov{\xi_t}+\xi_t)=2C_t\xi_t$
and $1-\rho_t=\xi_t(\ov{\xi_t}-\xi_t)=-2\zeta_4S_t\xi_t$. Let
$z=\omega'_0\xi_0\ov{\omega'_1\xi_1}$; then $|z|=1$.

For $i=0$ we get $\tau_{0,j}(D)=2\bigl(C_0\xi_0\omega'_0+(-1)^jC_1\xi_1\omega'_1\bigr)$, so
\[
4=|\tau_{0,j}(D)|^2=4\bigl(C_0^2+C_1^2+2(-1)^jC_0C_1\mathrm{Re}(z)\bigr)\qquad(j=0,1).
\]
Adding and subtracting the two equations gives $C_0^2+C_1^2=1$ and $C_0C_1\mathrm{Re}(z)=0$. For
$i=1$ we get $\tau_{1,j}(D)=-2\zeta_4\bigl(S_0\xi_0\omega'_0-(-1)^j\zeta_4S_1\xi_1\omega'_1\bigr)$,
and since $\mathrm{Re}(z\,\ov{\zeta_4})=\mathrm{Im}(z)$,
\[
4=|\tau_{1,j}(D)|^2=4\bigl(S_0^2+S_1^2-2(-1)^jS_0S_1\mathrm{Im}(z)\bigr)\qquad(j=0,1),
\]
whence $S_0^2+S_1^2=1$ and $S_0S_1\mathrm{Im}(z)=0$. As $z\ne0$, we have $C_0C_1=0$ or $S_0S_1=0$.

If $C_t=0$ for some $t$, then $C_{1-t}^2=1$ and hence $S_{1-t}=0$. If $S_t=0$ for some $t$, then
$S_{1-t}^2=1$ and hence $C_{1-t}=0$. So in all cases one of $\xi_0,\xi_1$ lies in $\{\pm1\}$ and the
other in $\{\pm\zeta_4\}$, i.e.\ $\{\rho_0,\rho_1\}=\{\xi_0^2,\xi_1^2\}=\{1,-1\}$.
\end{proof}

\begin{corollary}\label{cor:p2all}
Let $G$ be cyclic of order $4$ and let $h$ be even. Then every $\BH(G,h)$ matrix is of the form
\eqref{eq:family}.
\end{corollary}

\begin{proof}
Let $D\in\BH(G,h)$. By Proposition~\ref{prop:p2all},
$\{\rho_0,\rho_1\}=\{1,-1\}$ for the quotients $\rho_t=a_{g^{t+2}}/a_{g^t}$. Defining $\pi(t)\in\Z_2$
by $(-1)^{\pi(t)}=\rho_t$, we get $\pi(0)\ne\pi(1)$ and
$D=\sum_{t=0}^1a_{g^t}g^t\bigl(1+(-1)^{\pi(t)}g^2\bigr)$, which is \eqref{eq:family}.
\end{proof}

\begin{remark}\label{rem:p2all}
Corollary~\ref{cor:p2all} covers every even $h$, including $\nu_2(h)=1$.
Proposition~\ref{prop:p2all} uses only that the coefficients $a_x$ have absolute value $1$,
so its conclusion also applies to circulant complex Hadamard matrices of order $4$
whose entries need not be roots of unity.

Instead of the ad-hoc arguments given in this section, one may prove Corollary~\ref{cor:p2all} by applying the theory of vanishing sums to
obtain equations that the coefficients $a_x$ satisfy. That proof is more conceptual but we opt to give a more
elementary proof here.
\end{remark}

\section{Proof of the Main Theorem}\label{sec:proofmain}

\begin{lemma}\label{lem:ifpart}
Let $p$ be a prime with $p\mid h$, let $\omega_0,\dots,\omega_{p-1}\in\mu_h$, let $\pi$ be a
permutation of $\Z_p$, and let $D$ be given by \eqref{eq:family}. Then
\[
\tau_{i,j}(D)=p\,\omega_t\zeta_{p^2}^{-(i+pj)t}\qquad\text{with } t=\pi^{-1}(i)\qquad(i,j\in\Z_p).
\]
In particular $D$ is a $\BH(G,h)$ matrix.
\end{lemma}

\begin{proof}
The coefficients $\omega_t\zeta_p^{\pi(t)y}$ of $D$ lie in $\mu_h$, as $p\mid h$. We have
$\tau_{i,j}(g^{t+py})=\zeta_{p^2}^{-(i+pj)t}\zeta_p^{-iy}$, so
\[
\tau_{i,j}(D)=\sum_{t=0}^{p-1}\omega_t\zeta_{p^2}^{-(i+pj)t}\sum_{y=0}^{p-1}\zeta_p^{(\pi(t)-i)y}.
\]
The inner sum is $p$ if $\pi(t)\equiv i\pmod p$ and $0$ otherwise, and exactly one $t$ satisfies
$\pi(t)\equiv i$, namely $t=\pi^{-1}(i)$. Hence $|\tau_{i,j}(D)|^2=p^2$ for all $i,j$, and $D$ is
a $\BH(G,h)$ matrix.
\end{proof}

\begin{lemma}\label{lem:unique}
If $D=\sum_xa_xx$ is given by \eqref{eq:family}, then $\omega_t=a_{g^t}$ and
$\zeta_p^{\pi(t)}=a_{g^{t+p}}/a_{g^t}$ for all $t$. In particular, $D$ determines
$(\omega_0,\dots,\omega_{p-1})$ and $\pi$.
\end{lemma}

\begin{proof}
The coefficient of $g^{t+py}$ in \eqref{eq:family} is $\omega_t\zeta_p^{\pi(t)y}$; take $y=0$ and
$y=1$. Since $\pi(t)\in\Z_p$, it is determined by $\zeta_p^{\pi(t)}$.
\end{proof}

\begin{proof}[Proof of Theorem~\ref{thm:main}]
Let $D\in\BH(G,h)$. By Corollary~\ref{cor:pdividesh}, $p\mid h$.
If $p=2$, Corollary~\ref{cor:p2all} gives \eqref{eq:family} directly with parameters in $\mu_h$.
Now let $p$ be odd, put $c=\nu_p(h)\ge1$, and set
$\tilde h=ph$ if $c=1$ and $\tilde h=h$ if $c\ge2$. Then $D\in\BH(G,\tilde h)$ and
$\nu_p(\tilde h)\ge2$. By Theorem~\ref{thm:nu2odd} if $\nu_p(\tilde h)=2$, or
Theorem~\ref{thm:nu3} if $\nu_p(\tilde h)\ge3$, $D$ has the form \eqref{eq:family}
with parameters in $\mu_{\tilde h}$. By Lemma~\ref{lem:unique},
$\omega_t=a_{g^t}\in\mu_h$, so the parameters in fact lie in $\mu_h$.
Conversely, when $p\mid h$, Lemma~\ref{lem:ifpart} shows that every element of
\eqref{eq:family} belongs to $\BH(G,h)$. Lemma~\ref{lem:unique} gives uniqueness of the parameters.
\end{proof}

For $p\mid h$, Lemmas~\ref{lem:ifpart} and~\ref{lem:unique} show that the map
$\bigl((\omega_t)_{t},\pi\bigr)\mapsto D$ given by \eqref{eq:family} is a bijection from
$\mu_h^p\times\mathrm{Sym}(\Z_p)$ onto $\BH(G,h)$; hence $|\BH(G,h)|=p!\,h^p$.

\section*{AI Statement}
The main result (Theorem \ref{thm:main}) and its proof were obtained by the authors.
Claude Fable 5.1 and ChatGPT-6 Astra were used for proof checking and simplification,
as well as for typesetting.
The authors reviewed the AI-assisted suggestions and take full responsibility for the content of the paper.


\begin{thebibliography}{99}

\bibitem{Ber}
S.~D.~Berman, On the theory of group codes, \emph{Kibernetika} \textbf{3} (1967), no.~1, 31--39;
English translation in \emph{Cybernetics} \textbf{3} (1967), no.~1, 25--31.

\bibitem{Cha}
P.~Charpin, Une g\'en\'eralisation de la construction de Berman des codes de Reed et Muller
$p$-aires, \emph{Comm. Algebra} \textbf{16} (1988), no.~11, 2231--2246.

\bibitem{DS}
T.~Do Duc and B.~Schmidt, Bilinear forms on finite abelian groups and group invariant Butson
Hadamard matrices, \emph{J. Combin. Theory Ser. A} \textbf{166} (2019), 337--351.

\bibitem{DS26} T.~Do Duc and B.~Schmidt, Classification of Butson matrices invariant under elementary
abelian groups of order $p^2$, \emph{J. Algebraic Combin.} \textbf{63} (2026), Article 7.
\url{https://doi.org/10.1007/s10801-025-01494-4}.

\bibitem{GS}
G.~Hiranandani and J.-M.~Schlenker, Small circulant complex Hadamard matrices of Butson type,
\emph{European J. Combin.} \textbf{51} (2016), 306--314.

\bibitem{LM}
P.~Landrock and O.~Manz, Classical codes as ideals in group algebras, \emph{Des. Codes Cryptogr.}
\textbf{2} (1992), 273--285.

\bibitem{LCZ}
K.~H.~Leung, J.~H.~E.~Chue and M.~Zhao, On vanishing sums and cyclic Butson matrices,
\emph{Des. Codes Cryptogr.} \textbf{93} (2025), 5143--5157.

\bibitem{Sch99}
B.~Schmidt, Cyclotomic integers and finite geometry, \emph{J. Amer. Math. Soc.} \textbf{12}
(1999), 929--952.

\bibitem{Sch19}
B.~Schmidt, A survey of group invariant Butson matrices and their relation to generalized bent
functions and various other objects, in: K.-U.~Schmidt and A.~Winterhof (eds.), \emph{Combinatorics
and Finite Fields: Difference Sets, Polynomials, Pseudorandomness and Applications}, Radon Ser.
Comput. Appl. Math. \textbf{23}, De Gruyter, Berlin, 2019, pp.~241--251.

\bibitem{SWX21} B.~Schmidt, D.~Q.~Wong and Q.~Xiang, Constructions of Butson Hadamard matrices
invariant under abelian $p$-groups, J. Combin. Theory Ser. A 181 (2021), 105433.

\bibitem{Tao05}
T.~Tao, An uncertainty principle for cyclic groups of prime order, \emph{Math. Res. Lett.}
\textbf{12} (2005), 121--127.

\end{thebibliography}
\end{document}